\documentclass[12pt]{amsart} \usepackage[margin=1.2in]{geometry} \usepackage{graphicx,latexsym} \usepackage{amsfonts, amssymb, amsmath, amsthm}

\newtheorem{theorem}{Theorem}[section] \newtheorem{proposition}[theorem]{Proposition}  \newtheorem{corollary}[theorem]{Corollary}

\theoremstyle{definition} \newtheorem{definition} [theorem]{Definition}

\theoremstyle{remark} \newtheorem{remark}[theorem]{Remark} 

 \numberwithin{equation}{section}

\begin{document}

\title[DSTFT and directional regularity]{Directional short-time Fourier transform and directional regularity}
%\subtitle{Do you have a subtitle?\\ If so, write it here}

\author[S. Atanasova]{Sanja Atanasova} \address{Faculty of Electrical Engineering and Information Technologies\\ Ss. Cyril and Methodius
University\\Rugjer Boshkovik bb\\ 1000 Skopje, Macedonia} \email{ksanja@feit.ukim.edu.mk}

\author[S. Pilipovi\' c]{Stevan Pilipovi\' c}
\address{Department of Mathematics and Informatics,
University of Novi Sad, Trg Dositeja Obradovi\'{c}a 4, 21000 Novi Sad, Serbia}
\email{stevan.pilipovic@dmi.uns.ac.rs}
\author[K. Saneva]{Katerina Saneva } \address{Faculty of Electrical Engineering and Information Technologies\\ Ss. Cyril and Methodius
University\\Rugjer Boshkovik bb\\ 1000 Skopje, Macedonia} \email{saneva@feit.ukim.edu.mk}

 \keywords{}

\begin{abstract}
We give some new results related  to the  directional short-time Fourier transform (DSTFT) and extend them on the  spaces $\mathcal
K_{1}(\mathbb R^{n})$ and $\mathcal K_{1}({\mathbb R})\widehat{\otimes}\mathcal U(\mathbb C^n)$ and their duals. Then, we define multi-directional STFT and,
for tempered distributions,
 directional regular sets and their complements, directional wave fronts. Different windows with mild conditions on their support show the invariance of these notions
 related to window functions. Smoothness of  $f$ follows from the assumptions of the directional regularity in any direction.
\end{abstract}

\maketitle

\section{Introduction}
\label{intro}
In multidimensional time-frequency analysis, wave fronts are useful concept when analyzing where, how and why one distribution is singular, and when observing the direction in which the singularity occurs. Also, wave fronts are one of the crucial elements in the recent studies of the theory of distributions because of their ability to control the product of distributions.

The motivation of this paper is coming from  \cite{grafakos}, where Grafakos and
Sansing developed a theory that merges the Radon transform and time-frequency theory, and introduced the concept of directionally sensitive time-frequency analysis.
 Let $ g\in\mathcal S(\mathbb R)$  be a non-zero window function, $(u,x,\xi)\in \mathbb S^{n-1}\times\mathbb
 R\times\mathbb R^n$, where $\mathbb S^{n-1}$ is the unit sphere and $f\in L^{1}(\mathbb R^n)$. Then,
\begin{equation}\label{GRF}g_{u, x, \xi}(t)=e^{2\pi i
 \xi(u\cdot t-x)}g(u\cdot t-x),\,t\in\mathbb R^n,\end{equation}  called  Gabor ridge functions, can be viewed as time-frequency analysis elements in the Radon domain. By pairing the
function $f$ with $g_{u, x, \xi}$, Grafakos and Sansing provided an idea to localize information in time, frequency and direction defining a
directionally sensitive variant of the STFT.
They have shown  that it is not possible to obtain an exact reconstruction of a signal using the Gabor ridge functions (\cite[Thrm. 1]{grafakos}),
and therefore they have modified their class of functions to the weighted Gabor ridge functions (see \cite{grafakos} for details).  Their results for directionally sensitive time-frequency decompositions in $L^{2}(\mathbb R^n)$ based on Gabor systems in $L^{2}(\mathbb R)$ are generalized in \cite{Ole},  by showing similar results for  discrete and continuous frames.

Giv in \cite{GIV} introduced
another transform which is also a directionally sensitive variant of the STFT, letting
 \[{g }_{u,x,\xi }\left(t\right)= e^{2\pi i t\cdot \xi}g \left({t\cdot u
-x}\right),\  \  \ t\in {{\mathbb R}}^n.\]
 Using these functions he defined the directional
short-time Fourier transform (DSTFT) and proved several orthogonality results and reconstruction formulas for it \cite{GIV}.

The aim of this paper is twofold. In the first part (Section 3), we consider the DSTFT and its adjoint on the exponential type distributions, as an extension of the results of two  of us (cf. \cite{HA}) for tempered distributions to distributions of exponential type $\mathcal K_1'(\mathbb R^n)$. In this part we give another result concerning an inversion formula for DSTFT in comparison to \cite{GIV} and \cite{HA}, where an additional integral over the unit sphere $\mathbb S^{n-1}$ appears in several important formulae.

In the second part of the paper we  give an extension, introducing the multi-directional short-time Fourier transform (Section 4). Moreover, by a simple transformation of coordinates, we simplify our exposition considering directions of orthonormal basis $e_1,...,e_k$ of $\mathbb R^k$ in the framework of $\mathbb R^n$. In this way we present our main aim, namely, the analysis of the regularity properties of a signal $f(t), t\in\mathbb R^n$,
being a tempered distribution, through the  knowledge of the short-time Fourier transform in direction of selected coordinates. In other words, we  introduce and analyze  the directional wave fronts which can have applications in the time-frequency analysis.

\section{Preliminaries}\label{pre}
\subsection{Notation} The Fourier transform of a function $f\in L^1(\mathbb R^n)$ is defined as $\mathcal{F}(f)(\xi)=\widehat
{f}(\xi)=\int_{\mathbb R ^n} e^{-2\pi i x \cdot\xi}f(x)dx, \, \xi\in\mathbb R^n.$ The translation
and modulation operators are given by $T_{x}f(\: \cdot \:)=f( \: \cdot \:-x)$ and $M_{\xi}f(\: \cdot \:)=e^{2\pi i \xi \:\cdot\:}f(\:\cdot\:),$ $
x,\xi\in\mathbb R ^n,$ respectively. The operators  $M_{\xi} T_x$ and $T_x M_{\xi}$ are called time-frequency shifts. The notation $\langle
f,\varphi\rangle$ means dual pairing, whereas $(f,\varphi)$ stands for the $L^{2}$ inner product.

%%%%%%%%%%%%%%%%%%%%%%%%%%%%%%%%%%%%%%%%%%%%%%%%%%%%%%%%%%%%%%%%%%%%%%%%%%%%%%%%%%%%%%%%%%%%%%%%%%%%%%%%%%5555

\subsection{Spaces}\label{spaces} \noindent The Schwartz space of rapidly decreasing smooth functions and its dual, the space of
tempered distributions, are denoted by ${\mathcal S}({{\Bbb
R}}^n)$ and ${\mathcal S'}({\Bbb R}^n)$, respectively,
\cite{scwartz}. Recall \cite{hasumi} that  the space of exponentially rapidly decreasing smooth functions  ${\mathcal K}_1(\mathbb R ^n)$
is the space that consists of  $\varphi\in C^{\infty}(\mathbb{R}^{n})$ for which all the norms $$\rho_k(\varphi):=\sup_{t\in {\mathbb R  ^n},
\
|\alpha|\leq k}e^{k|t|}|\varphi^{(\alpha)}(t)|, \ \ \ k\in{\mathbb N}_0,$$ are finite. It is an FS-space and therefore Montel and reflexive. Moreover,
the space ${\mathcal K}_1(\mathbb R ^n)$ is nuclear. The dual space  ${\mathcal K}'_1({\mathbb R}^n)$ consists of all distributions of
the exponential form $f=\sum_{|\alpha|\leq l}(e^{s |\: \cdot\: |}f_\alpha)^{(\alpha)}$, where $f_\alpha\in L^{\infty}(\mathbb{R}^{n})$ \cite{hasumi}.
Next, recall \cite{hasumi} that  $\mathcal U(\mathbb C ^n)$ is the space of entire functions such that
$\varphi\in\mathcal{U}(\mathbb{C}^{n})$ if and only if $$ \theta_{k}(\varphi):=\sup_{z \in \Pi_{k}}(1+|z|^{2})^{k/2}|\varphi(z)|<\infty, \ \
\forall
k\in{\mathbb N}_{0}, $$ where $\Pi_{k}$ is the tube $\Pi_{k}= \mathbb{R}^{n}+i[-k,k]^{n}.$ The dual space $\mathcal{U}'(\mathbb{C}^{n})$, known as the
space of Silva tempered ultradistributions (see \cite{SS}, \cite{Morimoto},\cite{Z}), contains the space of analytic functionals.

As it turns out, the Fourier transform is a topological isomorphism from ${\mathcal K}_1(\mathbb R ^n)$ onto $\mathcal U(\mathbb C ^n)$, and extends
to
a topological isomorphism (with respect to strong topologies) $\mathcal{F}:{\mathcal K}'_1({\mathbb R}^n)\to \mathcal{U}'(\mathbb{C}^{n})$, \cite{hasumi}, \cite{Z}.

Next, we introduce the topological tensor product space  $\mathcal K_{1}({\mathbb
R})\widehat{\otimes}\mathcal U(\mathbb C^n)$  derived as the completion of the tensor product $\mathcal K_{1}({\mathbb
R})\otimes\mathcal U(\mathbb C^n)$  in the $\pi$- topology, same as the completion in the $\varepsilon$-topology \cite{treves}. The topology of $\mathcal K_{1}({\mathbb R})\widehat{\otimes} \mathcal U(\mathbb C^n)$ is given by the family of the norms
\[\rho^{l}_{k}(\Phi):=\sup_{(x,z)\in {\mathbb R} \times \Pi_{k}}e^{k|x|}(1+|z|^{2})
^{k/2}\left|\frac{\partial^{l}} {\partial x^{l}}\Phi(x, z)\right|, \ \ k, l\in{\mathbb N}_0.\]

Its dual
$(\mathcal K_{1}({\mathbb R})\widehat\otimes \mathcal U(\mathbb C^n) )'=\mathcal K'_{1}({\mathbb R})\widehat\otimes \mathcal U'(\mathbb C^n)$ will be used in our definition of the DSTFT of exponential distributions as it contains the range of this transform (cf. Subsection 3.2). If a measurable function $F$ satisfies
  \[\left|F\left(x,z \right)\right|\le C e^{s|x|}{{{\rm (1+}\left|z\right|{\rm )}}^{s}}, \ \ \ \left(x,z
\right) \in \mathbb
 R\times\mathbb C^n,\]
 for some $s,\, C>0$, then
we shall identify $F$ with an element of $\mathcal K'_{1}({\mathbb R})\widehat\otimes \mathcal U'(\mathbb C^n)$ via
\begin{equation}\label{*}
\left\langle F,\Phi \right\rangle :=\int_{\mathbb R }\int_{\mathbb R ^n}
F\left(x,\xi+i\eta \right)\Phi \left(x,\xi+i\eta\right)d\xi dx,
\end{equation}
$z=\xi+i\eta,\,\xi, \eta\in\mathbb R^n,\,\Phi \in \mathcal K_{1}({\mathbb R})\widehat\otimes \mathcal U(\mathbb C^n) $. \eqref{*} holds due to the Cauchy integral theorem.

%%%%%%%%%%%%%%%%%%%%%%%%%%%%%%%%%%%%%%%%%%%%%%%%%%%%%%%%%%%%%%%%%%%%%%%%%%%%
%%%%%%%%%%%%%%%%%%%%%%%%%%%%%%%%%%%%%%%%%%%%%%%%%%%%%%%%%%%%%%%%%%%%%%%%%%%%

\subsection{The short-time Fourier transform}\label{STFT} Let $f\in L^{2}({\Bbb R^n})$. Recall that the short-time Fourier transform (STFT) of  $f$ with respect to a window function $g\in L^{2} ({\Bbb R^n})$ is given by \begin{equation}
\label{STFT} V_{g} f(x,\xi ): \ =\langle f(t), \overline{M_{\xi } T_{x} g(t)}\rangle_t =\int _{{\Bbb R}^n} f(t)\overline{g(t-x)}e^{-2\pi
i\xi\cdot t} \ dt,\,\,x,\xi \in {\Bbb R^n}.
\end{equation}

The adjoint of $V_g$, over $L^2(\mathbb R^{2n})$, is given by  \begin{equation*} V_{g}^{*}F(t)=\iint_{\mathbb R^{2n}}F(x, \xi)g(t-x)e^{2\pi i\xi \cdot
t}dxd\xi.\end{equation*}
If $g\neq 0$ and $\psi\in L^{2}(\mathbb R ^n) $ is a synthesis window for $g$, that is, those one for which $(g,\psi ) \neq 0$, then  any $f\in L^2(\mathbb
R^n)$
can be recovered from its STFT via the inversion formula
 \begin{equation} \label{inverzna} f(t)=\frac{1}{( g
,\psi ) } \iint \nolimits _{{\mathbb R}^{2n} } V_{g} f(x,\xi ) M_{\xi } T_{x} \psi(t) dxd\xi. \end{equation}

Whenever the generalized inner product in \eqref{STFT} is well-defined, the definition of $V_{g}f$ can be viewed in a larger classes than
$L^{2}(\mathbb R^n ) $. It is easy to show that if $g\in\mathcal S(\mathbb R^n)\backslash \{0\}$ is fixed window, then $V_{g}:\mathcal S(\mathbb
R^n)\to\mathcal S(\mathbb R^{2n})$ and $V_{g}^{*}:\mathcal S(\mathbb R^{2n})\to\mathcal S(\mathbb R^n)$ are continuous mappings. We refer to \cite{gr01}, \cite{GZ2} for the basic STFT theory.

Moreover, in \cite{KPSV} authors have shown that if $g\in\mathcal K_{1}({\mathbb R}^n)\backslash \{0\}$ then $V_{g}:\mathcal K_{1}({\mathbb
R}^n)\to \mathcal K_{1}({\mathbb R}^n)\widehat{\otimes} \mathcal U(\mathbb C ^n)$ and $V_{g}^{*}:\mathcal K_{1}({\mathbb R}^n)\widehat{\otimes}
\mathcal U(\mathbb C ^n)\to \mathcal K_{1}({\mathbb R}^n)$ are continuous mappings.

One can define the STFT of a distribution $f\in \mathcal K'_1(\mathbb R^n)$ (resp. $\mathcal S'(\mathbb R^n)$) with respect to a window
$g\in\mathcal K_{1}({\mathbb R^n})$ (resp. $g\in\mathcal S(\mathbb R^n)$) as \begin{equation}\label{STFTonS} V_{g}f(x,\xi)=\langle f,
\overline{M_{\xi}T_{x}g}\rangle.\end{equation}

%%%%%%%%%%%%%%%%%%%%%%%%%%%%%%%%%%%%%%%%%%%%%%%%%%%%%%%%%%%%%%%%%%%%%%%%%%%%%%%%%%%5
%%%%%%%%%%%%%%%%%%%%%%%%%%%%%%%%%%%%%%%%%%%%%%%%%%%%%%%%%%%%%%%%%%%%%%%%

\subsection{The directional short-time Fourier transform}\label{DSTFT functions}

Let $u\in\mathbb S^{n-1}$. The {directional short-time Fourier transform} (DSTFT) of an integrable function
$f\in L^{1}(\mathbb R^n)$ (or $f\in \mathcal D'_{L^1}(\mathbb R^n)$) with respect to $g\in\mathcal S(\mathbb R)$
is given by
\begin{equation}\label{DSTFT} DS_{g,u}f(x,\xi):=\int_{\mathbb R^n}f(t)\overline{g(u\cdot t-x)}e^{-2\pi i t\cdot \xi}dt=\left\langle
f(t),{\overline{g }_{u,x,\xi}}(t)\right\rangle_t,\end{equation} where $(x,
\xi)\in \mathbb
 R\times\mathbb R^n$.

One can show, by the  use of results of Gr\" ochenig \cite{gr01} (we will demonstrate this in the proof of Proposition \ref{dd33} of Section \ref{s4}), that for a non-trivial $g\in \mathcal S(\mathbb R)$, with synthesis window  $\psi\in \mathcal S(\mathbb R)$ and
$f\in L^{1}(\mathbb R^n)$, the following reconstruction formula holds pointwisely,
\begin{equation}\label{invDSTFT} f(t)=\frac{1}{( g, \psi )}\int_{\mathbb R^n}\int_{\mathbb R}{DS_{g, u}f(x,\xi)\psi_{u,x,\xi}(t)dxd\xi}.\end{equation}

The reconstruction formula (\ref{invDSTFT}) allow us define an operator that maps functions on $\mathbb
R\times\mathbb R^{n}$ to functions on $\mathbb{R}^{n}$ as superposition of functions ${g }_{u,x,\xi}$. Given $g\in\mathcal{S}(\mathbb{R})$, we
introduce the \emph{directional synthesis operator} as
\begin{equation} \label{synthesis} DS_{g, u}^{*} \Phi(t):= \int_{\mathbb R^n }\int_{\mathbb R} \Phi(x,\xi){g }_{u,x,\xi}(
t)dxd\xi, \  \  \ t\in\mathbb{R}^{n}. \end{equation}

Thus, the relation (\ref{invDSTFT}) takes the form $(DS^{*}_{\psi, u}\circ DS_{g, u})f=( g, \psi ) f$.

The authors in \cite{HA} have discussed the problem of extending the definition of DSTFT on the space of tempered distributions.
Here, we study the DSTFT in the context of the space $\mathcal K'_{1}(\mathbb R^n)$ of distributions of exponential
type.

 If $f\in \mathcal{K}_{1}(\mathbb{R}^{n})$ and $g \in
\mathcal{K}_{1}(\mathbb{R})$, then we immediately get that
(\ref{DSTFT}) extends to a holomorphic function in the second variable. This means that $DS_{g, u}f(x,z)$ is entire in $z\in \mathbb{C}^{n}$. We write in the
sequel $z=\xi+i\eta$ with $\xi,\eta\in \mathbb{R}^{n}$. Note also that if $\Phi\in \mathcal K_{1}({\mathbb R})\widehat{\otimes} \mathcal U(\mathbb C^n)$ and $g\in \mathcal{K}_{1}(\mathbb{R})$,
then, using  the Cauchy theorem, we may write $DS^{\ast}_{g, u}\Phi$ as \begin{equation}\label{adjoint} DS_{g, u}^{*}\Phi(t) = \int_{\mathbb R
^{n}}\int_{\mathbb R}\Phi(x, \xi+i\eta)g(u\cdot t-x) e^{2\pi i(\xi+i\eta) \cdot t}dxd\xi, \end{equation}
 for arbitrary $\eta\in\mathbb{R}^{n}$.
In the next section, we will show that if $g \in
\mathcal{K}_{1}(\mathbb{R})$, then $ DS_{g, u}^{*}$ maps continuously $\mathcal K_{1}({\mathbb R})\widehat{\otimes} \mathcal U(\mathbb C^n)\to {{\mathcal K_{1}}}({{\mathbb R}}^n).$ It will then be shown that $ DS_{g, u}^{*}$ can be even extended to act on the distribution space $\mathcal K'_{1}({\mathbb R})\widehat\otimes \mathcal U'(\mathbb C^n)$.

As a simple consequence of Fubini's theorem, if $g\in\mathcal{K}_1(\mathbb{R})$, $f\in L^{1}(\mathbb{R}^{n})$ and $\Phi\in \mathcal K_{1}({\mathbb R})\widehat{\otimes} \mathcal U(\mathbb C^n)$,
 then one can easily prove \begin{equation} \label{transpose1} \int_{\mathbb{R}^{n}}f(t) DS^{*}_{g, u} \Phi(t)dt=\int_{\mathbb R^n}\int_{\mathbb R}
DS_{\overline{g}, u}f(x,\xi+i\eta)\Phi(x,\xi+i\eta)dxd\xi, \end{equation} $\xi, \eta\in\mathbb R^n$, and this can be written as
$ \left\langle f,DS_{\bar{g}, u}^{*}\Phi\right\rangle=\left\langle DS_{g, u} f,\Phi\right\rangle $ using (\ref{*}). As in \cite{HA}, we use this dual relation when defining the DSTFT of exponential distributions.

\section{The directional STFT of distributions of exponential type}
\subsection{Continuity of the directional STFT on ${{\mathcal K_{1}}}({{\mathbb R ^n}})$}\label{dir test functions}

Let $g\in\mathcal K_1(\mathbb R)\setminus\{0\}$. Then the DSTFT $DS_{g, u}$ is injective and $DS_{g, u}^{*}$ is surjective, due to the reconstruction formula (\ref{invDSTFT}).

Notice that we can extend the definition of the DSTFT as a sesquilinear mapping $ DS: (\varphi, g)\mapsto DS_{g, u}\varphi, $ $\varphi\in\mathcal K_1(\mathbb R^n)$, $g\in\mathcal K_1(\mathbb R)$, whereas the directional  synthesis
operator extends to the bilinear form $ DS^{*}: (\Phi, g)\mapsto DS^{*}_{g, u}\Phi, $ $\Phi\in {{\mathcal K_{1}}}({{\mathbb
R}})\hat{\otimes}{\mathcal U}(\mathbb C^{n})$.

\begin{theorem}\label{continuity theorem ridgelet}The mapping ${{ DS}}:{{\mathcal K_{1}}}({{\mathbb R}}^n)\times{{\mathcal K_{1}}}({{\mathbb
R}})\to {{\mathcal K_{1}}}({{\mathbb R}})\hat{\otimes}{\mathcal U}(\mathbb C^{n}) $ is continuous.
\end{theorem}

\proof

We will show that for given $k, l \in {{\mathbb N}}_0$, there exist $\nu,\tau \in \mathbb N_0$ and $C>0$ such that
\begin{equation}\label{GrindEQ__3_1_} \rho^{l}_{k}(DS_{g, u}\varphi)\leq C \rho_{\nu}(\varphi)\rho_{\tau}(g), \ \ \
\varphi\in\mathcal{K}_1(\mathbb{R}^{n}),\,g\in\mathcal{K}_1(\mathbb{R}). \end{equation}
Indeed, we have
\begin{align*}& e^{k|x|}(1+|z|^{2})^{k/2}\left|\frac{\partial^l}{\partial x^l}{DS}_{g, u}\varphi\left(x,z \right)\right|\\&
 =e^{k|x|}(1+|\xi+i\eta|^{2})^{k/2}\left|\frac{\partial^l}{\partial x^l}\int_{\mathbb
R^n}\varphi(t)\overline{g(u\cdot t-x)}e^{-2\pi i t\cdot (\xi+i\eta)}dt\right| \\ & \leq  e^{k|x|}(1+|\xi|^{2})^{k/2}(1+|\eta|^{2})^{k/2}\left|\int_{\mathbb
R^n}\varphi(t)\overline{g^{(l)}(u\cdot
t-x)}(-1)^{l}e^{-2\pi i t\cdot \xi}e^{2\pi t\cdot\eta}dt\right|\\&\leq
Ce^{k|x|}(1+nk^2)^{k/2}\left|\int_{\mathbb R^n}\varphi(t)\overline{g^{(l)}(u\cdot
t-x)}(-1)^{l}(1-\triangle _t)^{k/2}(e^{-2\pi i t\cdot \xi})e^{2\pi t\cdot\eta}dt\right|\\&\leq
Ce^{k|x|}(1+nk^2)^{k/2}\int_{\mathbb R^n}\left|(1-\triangle_t)^{k/2}\left(\varphi(t)\overline{g^{(l)}(u\cdot
t-x)}\right)\right|\left|e^{2\pi t\cdot\eta}\right|dt\\ & \leq Ce^{k|x|}(1+nk^2)^{k/2}\sum_{|k_1|+|k_2|=k}\left(
                                                                \begin{array}{c}
                                                                  k \\
                                                                  k_1, k_2 \\
                                                                \end{array}
                                                              \right)
\int_{\mathbb R^n}\left|\varphi^{(k_1)}(t)\overline{g^{(l+k_2)}(u\cdot t-x)}\right|e^{2\pi k|t|}dt\\&\leq
\widetilde{C_{k}}\sum_{|k_1|+|k_2|=k}\left(
                                                                \begin{array}{c}
                                                                  k \\
                                                                  k_1, k_2 \\
                                                                \end{array}
                                                              \right)
\int_{\mathbb R^n}e^{k|x-u\cdot t+u\cdot t|}\left|\varphi^{(k_1)}(t)\overline{g^{(l+k_2)}(u\cdot t-x)}\right|e^{2\pi
k|t|}dt\\&\leq \widetilde{C_{k}}\sum_{|k_1|+|k_2|=k}\left(
                                                                \begin{array}{c}
                                                                  k \\
                                                                  k_1, k_2 \\
                                                                \end{array}
                                                              \right)
\int_{\mathbb R^n}e^{k|u\cdot t-x|}e^{k|t|}\left|\varphi^{(k_1)}(t)\overline{g^{(l+k_2)}(u\cdot t-x)}\right|e^{2\pi k|t|}dt\\ &
= \widetilde{C_{k}}\sum_{|k_1|+|k_2|=k}\left(
                                                                \begin{array}{c}
                                                                  k \\
                                                                  k_1, k_2 \\
                                                                \end{array}
                                                              \right)
\int_{\mathbb R^n}e^{k|u\cdot t-x|}e^{(1+2\pi)k|t|}\left|\varphi^{(k_1)}(t)\right|\left|\overline{g^{(l+k_2)}(u\cdot t-x)}\right|dt.
\end{align*}
 %{\bf explain 8k-we added one more row at the end.}
\qed

%%%%%%%%%%%%%%%%%%%%%%%%%%%%%%%%%%%%%%%%%%%%%%%%%%%%%%%%%%%%%%%%%%%%%%%%%%5
%%%%%%%%%%%%%%%%%%%%%%%%%%%%%%%%%%%%%%%%%%%%%%%%%%%%%%%%%%%%%%%%%%%%%%%%%%%%%

We now analyze the directional synthesis operator.

\begin{theorem} \label{continuity theorem inverse} The bilinear mapping ${ DS}^{*}: ({{\mathcal
K_{1}}}({{\mathbb R}})\hat{\otimes}{\mathcal U}(\mathbb C^{n}))\times{{\mathcal K_{1}}}({{\mathbb R}})\to {{\mathcal K_{1}}}({{\mathbb R}}^n)$ is
continuous. \end{theorem}

\proof Let $g\in {{\mathcal K_{1}}}({{\mathbb R}})$, $\Phi\in {{\mathcal K_{1}}}({{\mathbb
R}})\hat{\otimes}{\mathcal U}(\mathbb C^{n})$ and $\varphi=DS_{g, u}^*\Phi $. Let  $\widehat{\Phi}_1(z_1,z)$ denote the Fourier
transform of $\Phi(x,z)$ with respect to the first variable and $\mathcal{F}^{-1}_{2}(\Phi)(x,t)$ denote the inverse Fourier
transform of $\Phi(x,z)$ with respect to the second  variable. We remark that $\widehat{\Phi}_1(z_1,z)$ is an entire function in $z_1=\omega+i\mu$, $\omega,\,\mu\in \mathbb R$. An application of the Cauchy theorem and the  Parseval formula give
$$\int_{\mathbb R} \Phi(x,z)g(u\cdot t-x) dx=\int_{\mathbb R}\widehat{\Phi}_1(\omega+i\mu,z)e^{-2\pi i (\omega+i\mu) u\cdot t}\widehat{g}(\omega+i\mu ) d\omega.$$

Observe that \begin{align*} &  \varphi(t)= DS_{g, u}^{*}\Phi(t) = \int_{\mathbb
R}\int_{\mathbb R^n} \Phi(x, \xi+i\eta)g(u\cdot t-x) e^{2\pi i(\xi+i\eta) \cdot t}d\xi d x \\&
=\int_{\mathbb
R}\int_{\mathbb R^n} \left(e^{-2\pi i (\omega+i\mu) u\cdot t}\widehat{\Phi}_1(\omega+i\mu,\xi+i\eta){\widehat{g}(\omega+i\mu)}\right)e^{2\pi i t\cdot
(\xi+i\eta)}d\xi d \omega\\&=
\int_{\mathbb R} e^{-2\pi i (\omega+i\mu) u\cdot t}{\widehat{g}(\omega+i\mu )}\left(\int_{\mathbb
R^n}\widehat{\Phi}_1(\omega+i\mu,\xi+i\eta)e^{2\pi i t\cdot (\xi+i\eta)}d\xi
\right)d \omega\end{align*}\begin{align*} &  =\int_{\mathbb R} e^{-2\pi i (\omega+i\mu) u\cdot t}{\widehat{g}(\omega+i\mu )}\mathcal{F}^{-1}_{2}{(\hat{\Phi}_{1})}(\omega+i\mu,t)d
\omega.
\end{align*}
 Hence, \begin{align} \label{Fourier} \nonumber  &\widehat{\varphi}(z)=\int_{\mathbb R}\widehat{g}(z_1)\widehat{\Phi}_1(z_1,z)*\widehat{e^{-2\pi i z_1 u\cdot t}}d\omega \\ & = \int_{\mathbb R}\widehat{g}(z_1)\widehat{\Phi}_1(z_1,z)*\delta({z_1}u+z)d\omega=\int_{\mathbb
R}\widehat{\Phi}_1(z_1,{z_1}
u+z)\widehat{g}(z_1)d\omega , \end{align} $z=\xi+i\eta\in\mathbb{C}^n, z_1=\omega+i\mu\in\mathbb C.$

We now prove the continuity of the bilinear directional synthesis mapping. Since the Fourier transform $g\mapsto\widehat{g}$ is a topological
isomorphism from $\mathcal K_{1}(\mathbb R^n)$ onto $\mathcal U(\mathbb C^{n})$, the family of seminorms $$\sigma_{k}(g)=
\theta_{k}(\widehat{g}), \  \ \  g\in\mathcal{K}_1(\mathbb{R}^n), \  \  \  k\in\mathbb N_{0}, $$ is a bases of seminorms for the topology of
$\mathcal{K}_1(\mathbb{R}^n)$.  %We will define a different family of seminorms on $\mathcal{K}_1(\mathbb{R}^n)$.

 %The seminorms $\dot{\rho}_{N}$,
%by $$ \dot{\rho}_{N}(\phi):=\sup_{(u,\omega)\in\mathbb{S}^{n-1}\times\mathbb{R}} \left|(1+|z|)^{N} \widehat{\phi}(z)\right|, \   \   \
%N\in\mathbb{N}_{0}, $$ is a base of continuous seminorms for the topology of $\mathcal{K}_1(\mathbb{R}^n)$.
We also know that the Fourier transform with respect to the first variable, $\Phi(x,z) \to \hat{\Phi}_1(z_1,z)$, is a topological isomorphism
from $\mathcal K_{1}(\mathbb R)\hat{\otimes}\mathcal U(\mathbb C^{n})$ onto $\mathcal
U(\mathbb C)\hat{\otimes}\mathcal U(\mathbb C^{n})$. Therefore, the family of seminorms $$\theta_{l,k}(\Phi)=\sup_{(z_1,z)\in {\Pi_l^1 \times \Pi_{k}}}(1+|z_1|^{2})^{l/2}(1+|z|^{2})
^{k/2}\left|\hat\Phi_1(z_1, z)\right|, \ \ l, k\in{\mathbb N}_0,\  \ \   $$
$\Pi_l^1=\mathbb R+i[-l,l],$ is a bases of seminorms for the topology of $\mathcal K_{1}(\mathbb R)\hat{\otimes}\mathcal U(\mathbb C^{n})$ .

We show that for a given $N\in\mathbb{N}_{0}$ there is $C>0$ such that $$ \sigma_{N}\left(DS^{*}_{g, u}\Phi\right)\leq C \sigma_{N+2}(g){\theta}_{0,N}({\Phi}).$$

Now, setting again $\varphi({x}):=DS_{g, u}^*\Phi (x)$ and using the expression (\ref{Fourier}),
we get \begin{align*} & (1+|z|)^{N}\left| \hat{\varphi}(z)\right|=(1+|z|)^{N}\left|\int_{\mathbb
R}\hat{g}(z_1)\hat{\Phi}_{1}(z_1,z_1 u+z)d\omega\right|\\ & \leq \int_{\mathbb
R}|\hat{g}(z_1)||\hat{\Phi}_{1}(z_1,z_1 u+z)|(1+|z+z_1 u|)^{N}(1+|z_1|)^{N}d\omega\\ & \leq
\sigma_{N+2}(g){\theta}_{0,N}({\Phi})\int_{\mathbb R}\frac{1}{(1+|z_1|)^2}d\omega, \end{align*} where $z_1=\omega+i\mu\in\mathbb C$.
\qed

%%%%%%%%%%%%%%%%%%%%%%%%%%%%%%%%%%%%%%%%%%%%%%%%%%%%%%%%%%%%%%%%%%%%%%%%%%%%%%%55
%%%%%%%%%%%%%%%%%%%%%%%%%%%%%%%%%%%%%%%%%%%%%%%%%%%%%%%%%%%%%%%%%

\subsection{Directional STFT on $\mathcal K'_{1}({\mathbb{R}^{n}})$}\label{ssds}

Let $ u\in\mathbb S^{n-1}$. The continuity results allow us to define the DSTFT of $f\in\mathcal K'_{1}(\mathbb R^n)$ with respect to $g\in \mathcal K_{1}(\mathbb R)$ as the
element $DS_{g,u}f\in{{\mathcal K'_{1}}}({{\mathbb R}})\hat{\otimes}{\mathcal U'}(\mathbb C^{n})$ whose
action on test functions is given by
\begin{equation}\label{def11}\langle DS_{g,u}{f}, \Phi\rangle:=\langle f, DS^{*}_{\overline{g},u}{\Phi}\rangle, \ \ \
\Phi \in{{\mathcal K_{1}}}({{\mathbb R}})\hat{\otimes}{\mathcal U}(\mathbb C^{n}).
\end{equation}

Then, the directional synthesis operator $DS^{*}_{g,u} :{{\mathcal K'_{1}}}({{\mathbb R}})\hat{\otimes}{\mathcal U'}(\mathbb C^{n})\to\mathcal{K}'_{1}({\mathbb{R}^{n}})$ can be defined
as
\begin{equation}\label{def21}\langle DS^{*}_{g,u}F, \varphi\rangle:=\langle F, DS_{\overline{g},u}{\varphi}\rangle, \  \  \ F\in{{\mathcal K'_{1}}}({{\mathbb R}})\hat{\otimes}{\mathcal U'}(\mathbb C^{n}),\ \varphi \in \mathcal K_{1}({\mathbb
R^{n}}).\end{equation}

We immediately obtain:
\begin{proposition} \label{pr1} Let $g\in\mathcal K_{1}(\mathbb R)$. The
directional short-time Fourier transform $DS_{g,u}:\mathcal{K}'_{1}({\mathbb{R}^{n}})\to {{\mathcal K'_{1}}}({{\mathbb
R}})\hat{\otimes}{\mathcal U'}(\mathbb C^{n})$ and the directional synthesis operator $DS^{*}_{g,u} :$  \break${{\mathcal K'_{1}}}({{\mathbb R}})\hat{\otimes}{\mathcal U'}(\mathbb C^{n})\to\mathcal{K}'_{1}({\mathbb{R}^{n}})$ are  continuous
linear maps. \end{proposition}

\subsection{Direct directional STFT on $\mathcal S'({\mathbb{R}^{n}})$}\label{dd100}
In the previous subsection, as well as in \cite{HA}, directional STFT on $\mathcal K'_{1}({\mathbb{R}^{n}})$ and on $\mathcal S'({\mathbb{R}^{n}})$ are defined  as transposed mappings.

We will consider in this subsection a direct definition of  directional STFT on $\mathcal S'({\mathbb{R}^{n}})$ as follows. Let $g\in \mathcal S(\mathbb R),$  $u\in\mathbb S^{n-1}$, and $x\in\mathbb R$. Then
\begin{equation}\label{sdstft}
\mathcal S'(\mathbb R^n)\ni f(t)\mapsto f(t)\overline{g(t\cdot u-x)}\in\mathcal S'(\mathbb R^n),
\end{equation}
\noindent and
$$\mathcal S'(\mathbb R^n)
\ni f(t)\overline{g(t\cdot u-x)}\mapsto \mathcal F(f(t)\overline{g(t\cdot u-x)}\,)(\xi)\in \mathcal S'(\mathbb R^n) $$
defines $DS_{g, u}f(x,\xi)$.
\begin{proposition}\label{dd44}
The direct definition of DSTFT and the one given via the transposed mapping coincide.
\end{proposition}
\proof
Let $f\in\mathcal S'(\mathbb R^n)$ and  $(f_k)$ be a sequence from $\mathcal S'(\mathbb R^n)$ which converges to $f$ in $\mathcal S'(\mathbb R^n).$ Since both definitions agree on $f_k$, for every $k$, the assertion follows by the continuity.
\qed

%%%%%%%%%%%%%%%%%%%%%%%%%%%%%%%%%%%%%%%%%%%%%%%%%%%%%%%%%%%%%%%%%%%%%%%%%%%5
%%%%%%%%%%%%%%%%%%%%%%%%%%%%%%%%%%%%%%%%%%%%%%%%%%%%%%%%%%%%%%%%%%%%%%%%%%%%%%%%%%%%%5

\section{Multi-directional STFT}\label{s4}
We will extend the directional STFT, introducing the $k$-directional STFT, $1\leq k\leq n.$ The case $k=1$ is explained in the previous part of the paper.

Note that the $k-$th tensor product completed in $\pi$- or $\varepsilon$- topology $\mathcal K_1(\mathbb R)\hat \otimes...\hat \otimes$ $\mathcal K_1(\mathbb R)$
equals to $\mathcal K_1(\mathbb R^k).$ The same holds for $\mathcal S(\mathbb R^k)$. Below we will use notations $(\mathcal K_1(\mathbb R))^k=\mathcal K_1(\mathbb R)\times...\times\mathcal K_1(\mathbb R)$ and $(\mathcal S(\mathbb R))^k=\mathcal S(\mathbb R)\times...\times\mathcal S(\mathbb R).$

Let $1\leq k\leq n$. Let $u^k=(u_1,...,u_k)$ where $u_i, i=1,...,k$ are independent vectors of $\mathbb S^{n-1}$, and $x^k=(x_1,...,x_k)\in\mathbb R^k$.
Let the nontrivial functions
$g_1,...,g_k$ belong to $\mathcal K_1(\mathbb R)$
(resp.,  $ \mathcal S(\mathbb R)$), $g^k=g_1\cdot...\cdot g_k\in (\mathcal K_1(\mathbb R))^k$ (resp.,  $(\mathcal S(\mathbb R))^k$) and $\xi\in\mathbb R^n.$

Let $f\in\mathcal K_1(\mathbb R^n)$ (resp., $\mathcal S(\mathbb R^n)$). Then, we define the $k$-directional STFT
by
\begin{equation}\label{dd2}
DS_{g^k,u^k} f(x^k,\xi): = \int_{\mathbb R^n} f(t)\overline{g_1(u_1\cdot t-x_1)}\cdot...\cdot \overline{g_k(u_k\cdot t-x_k)} e^{-2\pi it\cdot \xi}dt.
\end{equation}

\begin{proposition}\label{dd1111}
By (\ref{dd2}) is defined a continuous linear mapping of $\mathcal K_1(\mathbb R^n)$ (resp., $\mathcal S(\mathbb R^n)$) into
$ (\mathcal K_{1}(\mathbb R))^k\hat{\otimes}{\mathcal U}(\mathbb C^{n}) $
(resp., $ (\mathcal S(\mathbb R))^k\hat{\otimes}{\mathcal U}(\mathbb C^{n})$).

In particular, when $k=n$, (\ref{dd2}) is the short-time Fourier transform.
\end{proposition}

\proof
Let $A=[u_{i,j}]$ be a $k\times n$ matrix with raws $u_i, i=1,...,k$ and $I_{n-k,n-k}$ be the identity matrix. Let $B$ be an $n\times n$  matrix
determined by $A$ and $I_{n-k,n-k}$ so that $Bt=s$, where $$s_1=u_{1,1}t_1+...+u_{1,n}t_n,\; ... ,\; s_k=u_{k,1}t_1+...+u_{k,n}t_n,$$
$s_{k+1}=t_{k+1}, ..., s_n=t_n$. Clearly, it is regullar. Put $C=B^{-1}$ and $e^k=(e_1,...,e_k)$ where $e_1=(1,0,...,0),...,  e_k=(0,...,1)$ are
unit vectors of the coordinate system of $\mathbb R^k$.  Then, with the change of variables
$t=Cs$, and $\eta=C^t\xi$ ($C^t$ is the transposed matrix for $C$), one obtains, for $f\in\mathcal K_1(\mathbb R^n),$

$$DS_{g^k,u^k}f(x^k,\xi)=(DS_{g^k,e^k}(|C|f(C\cdot)))(x^k,\eta)$$
\begin{equation}\label{dd22}
=\int_{\mathbb R^n}f(s)\overline{g_1(s_1-x_1)}\cdot...\cdot \overline{g_k(s_k-x_k)}e^{-2\pi i s\cdot\eta}ds,
\end{equation}
where $|C|$ is the determinant of $C.$

Now, we immediately see the proof of the theorem since $\tilde f(s)=|C|f(Cs), s\in\mathbb R^n$ is an element of $\mathcal K_1(\mathbb R^n)$.
\qed

Put
$$g^k _{u^k,\,x^k,\,\xi}(t)=g_1(u_1\cdot t-x_1)\cdot...\cdot g_k(u_k\cdot t-x_k)e^{2\pi i t\cdot\xi}, t\in\mathbb R^n.
$$
Let for $g_i\in\mathcal K_1(\mathbb R)$ (resp.,  $g_i\in\mathcal S(\mathbb R)$), $\mathcal K_1(\mathbb R)$
(resp.,  $\mathcal S(\mathbb R)$)
be the synthesis window, $i=1,...,k$ and let
$$( g^k, \psi^k )={( g_1, \psi_1 )}\cdot...\cdot {( g_k, \psi_k )}\neq 0.
$$

We will prove the inversion formula  for the multi-directional STFT:
\begin{proposition}\label{dd33} Let $f\in\mathcal K_1(\mathbb R^n)$
 (resp.,  $\mathcal S(\mathbb R^n)$),
$g^k, \psi^k\in(\mathcal K_1(\mathbb R))^k$ (resp., $(\mathcal S(\mathbb R))^k$). Then
\begin{equation}\label{2invDSTFT}
f(t)=\frac{1}{( g^k, \psi^k )}\int_{\mathbb R^n}\int_{\mathbb R^k}{DS_{g^k,u^k}f(
x^k,\xi)\psi^k_{u^k,x^k,\xi}(t)dx^{k}d\xi}\end{equation}
pointwisely.
\end{proposition}
\proof
The proof is the same as for the short-time Fourier transform (see \cite{gr01}, Theorem 3.2.1 and Corollary 3.2.3).
We will use, after the change of variables the representation
(\ref{dd22}). Let $\tilde f_i(\cdot)=|C|f_i(C\cdot), i=1,2$. Actually, by the Parseval identity we have that for given $f_1,f_2\in L^2(\mathbb R^n)$ and $g^k,\psi^k\in (\mathcal K_1(\mathbb R))^k$,
$$
 (DS_{g^k,u^k}f_1(x^k,\xi),DS_{\psi^k,u^k}f_2(x^k,\xi))_{L^2(\mathbb R^k \times \mathbb R^n)}
$$
\begin{equation}\label{ddd12}
=(DS_{g^k,e^k}\tilde f_1(x^k,\xi),DS_{\psi^k,e^k}\tilde f_2(x^k,\xi))_{L^2(\mathbb R^k \times \mathbb R^n)}
  =( \tilde f_1,\tilde f_2)_{L^2(\mathbb R^n)}( \overline{g^k},\overline{\psi^k})_{L^2(\mathbb R^k)}.
\end{equation}

We obtain the reconstruction formula
 (\ref{2invDSTFT}) as a consequence of (\ref{ddd12}), as in the quoted corollary of \cite{gr01}.
\qed

 Let $f\in \mathcal{K}_{1}(\mathbb{R}^{n})$. We have that
(\ref{dd2}) extends to a holomorphic function, i.e. \allowbreak $DS_{g^k,u^k}f(x^k,z)$ is entire in $z\in \mathbb{C}^{n}$. As in the case $k=1,$
 if $\Phi\in \mathcal
(\mathcal K_{1}({\mathbb R}))^k\widehat{\otimes} \mathcal U(\mathbb C^n)$
and $g^k\in (\mathcal{K}_{1}(\mathbb{R}))^k$,
 for arbitrary $\eta\in\mathbb{R}^{n}$ and  the Cauchy theorem, we can write
$$ DS_{g^k,u^k}^{*}\Phi(t) = \int_{\mathbb R
^{n}}\int_{\mathbb R^k}\Phi(x^k, \xi+i\eta)
$$
\begin{equation}\label{2adjoint}
g_1(u_1\cdot t-x_1)\cdot ...\cdot
g_k(u_k\cdot t-x_k) e^{2\pi i(\xi+i\eta) \cdot t}dx^{k}d\xi,\; \; t\in\mathbb R^n.
\end{equation}

Multi-directional STFT on dual spaces $\mathcal K_1'(\mathbb R^n)$ and $\mathcal S'(\mathbb R^n)$ can be defined as in the case $k=1$ (cf. Subsection 3.2 and 3.3).

The next theorem  connects multi-directional STFT's with respect to different windows. It is crucial for the main theorem of Section 5.

\begin{theorem}\label{d444} Let $u_1,...,u_k\in\mathbb S^{n-1}$ be independent.
Let  $h_1,...,h_k,$ $ g_1,..., g_k,$ $\gamma_1,...,\gamma_k$ belong to $\mathcal S(\mathbb R)$ where $\gamma_i$ is synthesis
window for $g_i, i=1,...,k. $ Let $f\in\mathcal S'(\mathbb R^n)$. Then
$$DS_{h^k,u^k}f(y^k,\eta)=(DS_{g^k,u^k}f(s^k,\zeta))*(DS_{h^k,u^k}\gamma^k(s^k,\zeta))(y^k,\eta).
$$
\end{theorem}
\proof By the use of the change of variables given in the  proof of Proposition \ref{dd1111}, it follows that it is enough to prove the assertion for $u^k=e^k.$ Let $F\in{{\mathcal S}'}({{\mathbb R^k}})\hat{\otimes}{\mathcal U'}(\mathbb C^{n})$.
Then
\begin{eqnarray*}&& DS_{h^k,e^k}(DS_{\gamma^k,e^k}^{*}F)(y^k,\eta) = \int_{\mathbb R^n}(\int_{\mathbb R
^{n}}\int_{\mathbb R^k}F(x^k, \xi)
\\&&
\gamma_1(t_1-x_1)\cdot...\cdot
\gamma_k(t_k-x_k) e^{2\pi i\xi \cdot t}dx^kd\xi)
\overline {h_1(t_1-y_1)}\cdot...\cdot\overline {h_k(t_k-y_k)} e^{-2\pi it\cdot\eta}dt
\\&=& \int_{\mathbb R^n}\int_{\mathbb R
^{k}}(\int_{\mathbb R^n}
\overline{h_1(t_1-(y_1-x_1))}\cdot...\cdot
\overline{h_k(t_k-(y_k-x_k))}
\\&&\gamma_1(t_1)\cdot...\cdot \gamma_k(t_k) e^{-2\pi it\cdot(\eta-\xi)}dt)F(x^k,\xi)dx^kd\xi\\
&=& \int_{\mathbb R^n}\int_{\mathbb R
^{k}}(\int_{\mathbb R^n}\gamma^k(t^k)
\overline{h^k(t^k-(y^k-x^k))}
 e^{-2\pi it\cdot(\eta-\xi)}dt)F(x^k,\xi)dx^kd\xi\\&=&\int_{\mathbb R^n}\int_{\mathbb R^{k}}F(x^k,\xi)DS_{h^k,e^k}\gamma^k(y^k-x^k,\eta-\xi)dx^kd\xi.
\end{eqnarray*}
Now, we put $F=DS_{g^k,e^k}f$ and obtain
\begin{equation}\label{dd333}
DS_{h^k,e^k}f(y^k,\eta)=(DS_{g^k,e^k}f(s^k,\zeta))*(DS_{h^k,e^k}\gamma^k(s^k,\zeta))(y^k,\eta).
\end{equation}
This completes the proof of the theorem.
\qed

%%%%%%%%%%%%%%%%%%%%%%%%%%%%%%%%%%%%%%%%%%%%%%%%%%%%%%%%%%%%%%%%%%%%%%%%%
%%%%%%%%%%%%%%%%%%%%%%%%%%%%%%%%%%%%%%%%%%%%%%%%%%%%%%%%%%%%%55
\section{Directional wave fronts}

DSTFT can be used in
the detection of singularities determined by the hyperplanes orthogonal to vectors $u_1,...,u_k.$ For this purpose, we introduce
(multi)-directional regular sets and wave front sets for tempered distributions using the direct multi-directional STFT.

The proofs of Proposition \ref{dd1111} and Theorem \ref{d444} show that we can simplify our exposition by the use of the linear transformation $C$ of Proposition \ref{dd1111}
and transfer the $u^k$- DSTFT to $e^k$-
DSTFT. Thus, in order to simplify our exposition of this section, we will consider regularity properties in the framework of the direction $u^k=e^k.$

If $k=1$, we consider direction $e^1=e_1$ while for $1<k\leq n$, we consider direction $e^k=(e_1,...,e_k).$
Let $k=1$ and $x_0=x_{0,1}\in\mathbb R$. Put
$
\Pi_{e^1,x_0,\varepsilon}=\Pi_{x_0,\varepsilon} :=\{t\in \mathbb R^n; |t_1- x_0|<\varepsilon\}.
$
It is a part of $\mathbb R^n$ between two hyperplanes orthogonal to $e_1$, that is,
$$
 \Pi_{x_0,\varepsilon}= \bigcup_{x\in (x_0-\varepsilon,x_0+\varepsilon)} P_{x},
\;\; \; (x_0=(x_0, 0,...,0)),$$
 and
$P_{x}$ denotes the hyperplane orthogonal to $e_1$ passing through  $x$.

We keep the notation of Section 4.
Put
$$\Pi_{e^k,x^k,\varepsilon}=\Pi_{e_1,x_1,\varepsilon}\cap...\cap\Pi_{e_k,x_k,\varepsilon}, \; \Pi_{e^k,x^k}=\Pi_{e_1,x_1}\cap...\cap\Pi_{e_k,x_k}.
$$
The first set is a paralelopiped determined by $2k$ finite edges while the other edges are infinite. The set  $\Pi_{e^k,x^k}$ equals  $\mathbb R^{n-k}$ translated by vectors $\vec{x_1},...,\vec{x_k}.$ We will call it $n-k$-dimensional element of  $\mathbb R^n $ and denote it as $P_{e^k,x^k}\in\mathbb R^k.$
 If $k=n$, then this is just the point $x^n=(x_1,...,x_n).$

\begin{definition}\label{wp}
Let  $f\in \mathcal S'(\mathbb R^n)$. It is said that $f$ is $k$-directionally  regular  at $(P_{e^k,x_0^k},\xi_0)\in
\mathbb R^n\times \mathbb R^n\setminus \{0\}$
if there exists $g^k\in (\mathcal D(\mathbb R))^k$, $g^k(0)\neq 0$, a product of
open balls $L_r(x^k_0)=L_r(x_{0,1})\times...\times L_r(x_{0,k})\in\mathbb R^k$
and a  cone $\Gamma_{\xi_0}$ such that  for every
$N\in\mathbb N$  there exists $C_N>0$
such that
\begin{equation*}\label{rhh}
\sup_{x^k\in L_r(x_0^k),\,\xi \in\Gamma_{\xi_0}}|DS_{g^k, e^k}f(x^k,\xi)|\end{equation*}\begin{equation}\label{rh}
=\sup_{x^k\in L_r(x_0^k),\,\xi \in\Gamma_{\xi_0}}|\mathcal F
(f(t)\overline{g^k(t^k-x^k)})(\xi)|\leq C_N(1+|\xi|^2)^{-N/2}.
\end{equation}
\end{definition}

Note that for $k=n$
we have classical H\" ormander's regularity.

\begin{remark}\label{wf}
a) If $f$ is $k$-directionally regular at  $(P_{e^k,x_0^k},\xi_0)$, then there exists an open ball (with radius $r$ and center $x_0^k$) $ L_r(x_0^k)$ and an open cone $\Gamma\subset\Gamma _{\xi_0}$
so that $f$ is $k$-directionally regular at  $(P_{e^k,z_0^k},\theta_0)$ for any $z_0^k\in L_{r}(x_0^k)$ and $\theta_0 \in \Gamma.$ This implies that the union  of all $k$-directional regular points $(P_{e^k,z_0^k},\theta_0)$, $(z_0^k,\theta_0)\in L_{r}(x_0^k)\times\Gamma$ is an open set of $\mathbb R^n\times\mathbb R^n\setminus\{0\}$.

b) Denote by $Pr_{k}$ the projection of $\mathbb R^n$ onto $\mathbb R^k$. Then, the  $k$-directional regular
point $(P_{e^k,x_0^k},\xi_0)$, considered in $\mathbb R^n\times\mathbb R^n\setminus\{0\}$ with respect to the first $k$ variables, equals $(Pr_k^{-1}\times I_\xi)(P_{e^k,x_0^k},\xi_0)$ ($I_\xi $ is the identity matrix on $\mathbb R^n$).

We define the $k$-directional wave front  as the complement in
$\mathbb R^k\times\mathbb R^n\setminus\{0\}$ of all    $k$-directional regular points $(P_{e^k,x_0^k},\xi_0)$. This set is denoted as
$WF_{e^k}f.$ In $\mathbb R^n\times \mathbb R^n\setminus \{0\},$ this is $(Pr_k^{-1}\times I_\xi)(WF_{e^k}f).$
\end{remark}

\begin{proposition}
Set
$WF_{e^k}(f)$
is closed  in $\mathbb R^k\times\mathbb R^n\setminus\{0\}$ (and $\mathbb R^n\times \mathbb R^n\setminus \{0\}$).
\end{proposition}

We will use notation $B_r(0^k)$ to denote a closed ball in $\mathbb R^k$ with center at  zero $0^k$  and radius $r>0.$
Our main theorem relates directional regular sets for two DSTFT's.

\begin{theorem} \label{nwr}  If (\ref{rh}) holds for some $g^k\in(\mathcal D(\mathbb R))^k$, then it holds for every $h^k\in(\mathcal D(\mathbb R))^k, $ $(h^k(0)\neq 0)$ supported by a ball
$B_\rho(0^k)$, where $\rho\leq\rho_0$ and $\rho_0$ depends on $r$ in  (\ref{rh}).
\end{theorem}
\proof

Since $g^k$ and $h^k$ are compactly supported,  the integration with respect to $x^k$, which will be performed below, is finite. Moreover, we can assume that $f$ is a
continuous polynomially bounded function. If not, let $f=P(D)F$, where $F$ is polynomially bounded and continuous while $P(D)$ is a differential operator with constant coefficients. In this case  we can perform partial integration and transfer the differentiation from $f$ on other factors of the integrand which do not effects the proof. So the analysis can be continued with $f$ continuous and polynomially bounded.

We use Proposition \ref{d444}, that is, the form (\ref{dd333}).  Assume that (\ref{rh}) holds and that $\gamma^k$ is chosen so that
 $\mbox{ supp }\gamma^k\subset B_{\rho_1}(0^k)$ and  $\rho_1<r-r_0$.
Let  $h^k\in(\mathcal D(\mathbb R))^k$ and $\mbox{ supp }h^k\subset B_{\rho}(0^k)$.
  We will find $\rho_0$ such that  (\ref{rh}) holds for $DS_{h^k,e^k}f(y^k,\eta)$, with
$y^k\in B_{r_0}(x_0^k), $ $\eta\in\Gamma_1\subset\subset \Gamma_{\xi_0},$
for $ \rho\leq \rho_0$
($\Gamma_1\subset\subset \Gamma_{\xi_0}$ means that $\Gamma_1\cap \mathbb S^{n-1}$
is a compact subset of  $\Gamma_{\xi_0}\cap \mathbb S^{n-1}$).

We need the next simple observations:
$$
|p^k|\leq \rho_1,\;\; |y^k-x_0^k|\leq r_0 \mbox{ and }\;\; |p^k-((y^k-x_0^k)-(x^k-x_0^k))|\leq \rho
$$
\begin{equation}\label{suporti}
\Rightarrow |x^k-x_0^k|\leq \rho+\rho_1+r_0.
\end{equation}
So, we choose $\rho_0$ such that
$\rho_0+\rho_1<r-r_0.$
Then
\begin{equation}\label{suporti2}\rho+\rho_1+r_0<r \;\mbox{ holds for }\;
 \rho\leq\rho_0.
\end{equation}

Let $\Gamma_1\subset\subset \Gamma_{\xi_0}$.
Then, with a suitable $c\in (0,1)$,
\begin{equation}\label{gam}
\eta\in \Gamma_1, |\eta|>1 \mbox{ and } |\eta-\xi|\leq c|\eta|\Rightarrow \xi \in\Gamma_{\xi_0}; \;\;
|\eta-\xi|\leq c|\eta|\Rightarrow |\eta|\leq (1-c)^{-1}|\xi|.
\end{equation}
Let $y^k\in B_{r_0}(x^k_0), \eta \in\Gamma_1$.
Integrals which will appear below are considered as oscillatory integrals. We have
$$
|DS_{h^k,e^k}f(y^k,\eta)|
=\left|\int_{\mathbb R^k}\int_{\mathbb R^n}DS_{g^k,e^k}f(x^k,\xi)DS_{h^k,e^k}\gamma^k(y^k-x^k,\eta-\xi)d\xi dx^k\right|.
$$
Consider $$ J_1=\int_{\mathbb R^n}DS_{g^k,e^k}f(x^k,\eta-\xi)d\xi \mbox{ and }\;  J_2=\int_{\mathbb R^n}DS_{h^k,e^k}\gamma^k(y^k-x^k,\xi)d\xi.$$
Then,  by the use of partial integration, we have, in the oscillatory sense, (with the assumption that $n $ is odd),
$$J_1=\int_{\mathbb R^n}\int_{\mathbb R^n}
\frac{f(t)}{(1+|2\pi t|^2)^{(n+1)/2}}\overline{g^k(t^k-x^k)}(1-\Delta_\xi)^{(n+1)/2}e^{-2\pi i t\cdot(\eta-\xi)}dtd\xi
$$
(for $n$ even, we take $n+2$ instead of $n+1$). This integral still diverges with respect to $\xi$, while $J_2$ converges, since
\begin{eqnarray*}J_2&=&\int_{\mathbb R^n}\int_{B_{\rho_1}(0^k)}
\frac{\gamma^k(p^k)\overline{h^k(p^k-(y^k-x^k))}}{(1+|2\pi \xi|^2)^{s/2}}(1-\Delta_p)^{s/2}e^{-2\pi i p\cdot\xi}dpd\xi
\\&=&
\int_{\mathbb R^n}\int_{B_{\rho_1}(0^k)}
(1-\Delta_p)^{s/2}\frac{\gamma^k(p^k)\overline{h^k(p^k-(y^k-x^k))}}{(1+|2\pi \xi|^2)^{s/2}}e^{-2\pi i p\cdot\xi}dpd\xi.
\end{eqnarray*}
Rewrite
$$|DS_{h^k,e^k}f(y^k,\eta)|=\int_{\mathbb R^k}|(\int_{|\eta-\xi|\leq c|\eta|}+\int_{|\eta-\xi|\geq c|\eta|})(...)d\xi|dx^k=I_1+I_2.
$$
Then,
$$I_1\leq \int_{\mathbb R^k}(\sup_{|\eta-\xi|\leq c|\eta|}
|DS_{g^k,e^k}f(x^k,\eta-\xi)\int_{|\eta-\xi|\leq c|\eta|}
|DS_{h^k,e^k}\gamma(y^k-x^k,\xi)|d\xi)dx^k.
$$
Now, we use (\ref{suporti}), (\ref{suporti2}) and $(1+|\eta|^2)^{N/2}\leq C(1+|\xi|^2)^{N/2}$, for $|\xi|\geq (1-c)|\eta|$. This implies
$$\sup_{y^k\in B_{r_0}(x^k_0),\,\eta\in\Gamma_1}(1+|\eta|^2)^{N/2}I_1\leq \int_{B_r(x_0^k)}\left (\sup_{\xi\in \Gamma_{\xi_0}}
|DS_{g^k,e^k}f(x^k,\xi)|(1+|\xi|^2)^{N/2}\right .
$$
$$\left .\times \int_{|\xi|\geq (1-c)|\eta|}
|DS_{h^k,e^k}\gamma^k(y^k-x^k,\xi)|d\xi\right )dx^k.
$$
Now by the finiteness of $J_2$, we obtain that $I_1$ satisfies the necessary estimate of (\ref{rh}). Let us consider $I_2.$
$$I_2\leq \int_{\mathbb R^k}\left|\int_{|\xi|\geq c|\eta|}
DS_{g^k,e^k}f(x^k,\eta-\xi)DS_{h^k,e^k}\gamma^k(y^k-x^k,\xi)d\xi\right| dx^k.
$$

Let $K=\{\xi: |\xi|\geq c|\eta|\}$.
Denote by $\widetilde \kappa_d, 0<d<1,$ the characteristic function of
 $K_{d}=\bigcup_{\xi\in K}L_d(\xi)$, that is, $K_d$ is open $d$-neighborhood of $K.$
Then, put $$\kappa_\eta=\widetilde{\kappa}_{d}*\varphi_{d},$$
where $\varphi_d=\frac{1}{d^n}\varphi(\cdot/d)$,
$\varphi\in\mathcal D(\mathbb R^n)$ is non-negative, supported by the ball $B_1(0)$ and equals
$1/2$ on $B_{1/2}(0).$ This construction implies that $\kappa_\eta$ equals one on $K$,  is supported by $K_{2d}.$
Moreover,  all the derivatives of $\kappa_\eta$ are bounded.    Assume that $n$ is odd and that $s$ is even and sufficiently large.
Then,
$$\sup_{y^k\in B_{r/2}(x^k_0),\,
 \eta\in \Gamma_1}I_2\leq C \int_{\mathbb R^k}|\int_{\mathbb R^n}\kappa_\eta(\xi)
DS_{g^k,e^k}f(x^k,\eta-\xi)DS_{h^k,e^k}\gamma^k(y^k-x^k,\xi)d\xi|dx^k
$$
\noindent$$\leq C\int_{\mathbb R^k} |\int_{\mathbb R^n_\xi}(\int_{\mathbb R^n_t}\frac{f(t)}{(1+|2\pi t|^2)^{(n+1)/2}}\overline{g^k(t^k-x^k)}e^{-2\pi i t\cdot(\eta-\xi)}dt)(1-\Delta_\xi)^{(n+1)/2}
$$
$$\big(\frac{\kappa_\eta(\xi)}{(1+|2\pi \xi|^2)^{s/2}} \int_{\mathbb R^n_p}\gamma^k(p^k)\overline{h^k(p^k-(y^k-x^k))}(1-\Delta_p)^{s/2} e^{-2\pi ip\cdot\xi}dp\big) d\xi|dx^k.
$$
Then, for every $y^k\in B_{r_0}(x_0^k)$ and $\eta\in\Gamma_1$, choosing $s>N+n$ (and being even) as well as using the Petree
inequality, we see that all the integrals on the right hand side of
$$(1+|\eta|^2)^{N/2}I_2\leq C\int_{\mathbb R^k}\int_{\mathbb R^n_\xi}
(\int_{\mathbb R^n_t}\frac{|f(t)|}{(1+|2\pi t|^2)^{(n+1)/2}}|\overline{g^k(t^k-x^k)}|dt
$$
$$
\big(\frac{(1+|\xi|^2)^{N/2}}{(1+|\eta-\xi|^2)^{N/2}}(1-\Delta_\xi)^{(n+1)/2}
(\frac{\kappa_\eta(\xi)}{(1+|2\pi \xi|^2)^{s/2}}\big)
$$$$ \int_{\mathbb R^n_p}|(1-\Delta_p)^{s/2}\big(\gamma^k(p^k)\overline{h^k(p^k-(y^k-x^k))}\big)|dp \big)d\xi dx^k
$$
are finite. This completes the proof of the theorem.
\qed
\begin{remark}\label{ob}
If $\mbox{supp }g^k\subset B_{a}(0^k)$, then we see that (\ref{rh}) shows the behaviour of
$f(t), t\in Pr_k^{-1}(B_{a+r}(x_0^k))$ in the direction of $\xi_0.$
\end{remark}
\begin{corollary}
\label{suz}
Let $g^k\in(\mathcal S(\mathbb R))^k$, supported by a ball $B_a(0^k)$, have  synthesis window $\gamma^k$ supported by $B_{\rho_1}(0^k),$
$\rho_1\leq a.$ Assume that in (\ref{rh}) we have $2r$ instead of $r$, that is,
\begin{equation}\label{2rh}
\sup_{x^k\in L_{2r}(x_0^k),\,\xi \in\Gamma_{\xi_0}}|DS_{g^k, e^k}f(x^k,\xi)|\leq C_N(1+|\xi|^2)^{-N/2}.
\end{equation}
 Moreover, assume that $a<r.$ Then, for any $h^k\in\mathcal (D(\mathbb R))^k$ with support $B_{\rho}(0^k)$, $\rho<a$, there exists $r_0$ and $\Gamma_1\subset\subset \Gamma_{\xi_0}$ such that (\ref{2rh}) holds for $DS_{h^k,e^k}f(y^k,\eta)$ with the supremum over
$y^k\in B_{r_0}(x_0^k)$ and $ \eta \in \Gamma_1$.
\end{corollary}
\proof
With the notation of
Theorem \ref{nwr}, we have, similarly
as in (\ref{suporti}),
$$|x^k-x_0^k|\leq \rho+\rho_1+r_0<\rho+a-r_0+r_0=a+\rho<2r.
$$
This implies $|x^k-x^k_0|<2r,$ so that the supremum in the estimate of $I_1$ holds. In the same way as in Theorem \ref{nwr}, we perform the proof.
\qed

\begin{theorem}\label{w1}
If
$(P_{e^k,x^k},\xi)$ is a $k$-directional regular point of $f\in\mathcal S'(\mathbb R^n)$ for every $\xi\in\mathbb R^n\setminus \{0\}$,
then  $ f\in\mathcal E(\mathbb R^n).$
\end{theorem}
\proof
Assume that %$g\in\mathcal D(\mathbb R)$
 supp $g^k \subset B_\rho(0^k).$
Let $x_0^k\in \mathbb R^k.$
For every $\xi\in \mathbb S^{n-1}$ there exist a ball
$L^\xi_{r}(x_0^k)$ and a cone $\Gamma_\xi$ such that (\ref{rh}) holds. As in the classical theory, the compactness of $\mathbb S^{n-1}$
 implies that there exists $r>0$ such that for every $N>0$ there exists $C_N>0$ such that% for the direction $u_0\in\mathbb S^{n-1},$
$$\sup_{x^k\in L_{r}(x_0^k)}|\mathcal F(f(t)\overline{g_1(t_1-x_1)}\cdot...\cdot\overline{g_k(t_k -x_k)})(\xi)|\leq C_N(1+|\xi|^2)^{-N/2},\xi\in\mathbb R^n.$$
Thus, $f(t)\overline{g^k(t^k-x^k)}\in\mathcal E(\mathbb R^n)$ for every $x^k\in L_r(x_0^k)$. Since $|t^k-x^k|< \rho$, we see that
$t$ must lie in some $Pr^{-1}_k(L_{r+\rho}(x_0^k)).$ Thus, for every point of $\mathbb R^n$ there exists an open set around it where $f$ is smooth.
This completes the proof of the theorem.
\qed

%\begin{acknowledgements}
%S. Atanasova and K. Saneva gratefully acknowledge support by Macedonian-Austrian bilateral project through the grant 10-1491/2.
%\end{acknowledgements}

% Non-BibTeX users please use
 
\end{document}